\documentstyle[12pt]{article}
\def \reel{ {\rm I}\!{\rm R} }
 \newcommand{\too}{\longrightarrow}
\newcommand{\om}{\omega}

\newcommand{\na}{\nabla}
\newcommand{\wi}{\widetilde}
\newcommand{\al}{\alpha}
\newcommand{\be}{\beta}
\newcommand{\ga}{\gamma}

 \def \rat{ {\rm Q}\kern-.65em {}^{{}_/ }}

\newtheorem{Def}{D\'efinition}[section]

\newtheorem{pr}{Proposition}[section]

\title{Compatibilit\'e des structures pseudo-riemanniennes et des structures de
Poisson} \author{M. Boucetta}
\date{}
\begin{document}	
\maketitle		
{\bf R\'esum\'e.} Dans cette note, nous allons introduire deux notions naturelles et ind\'ependantes de compatibilit\'e entre une m\'etrique pseudo-riemannienne et une structure de Poisson et ce en utilisant la notion de d\'eriv\'ee contravariante introduite dans [1], nous allons donner les premi\`eres propri\'et\'es de couple de m\'etrique et de structure de Poisson compatible pour l'une o\`u l'autre  et finalement, nous allons donner des exemples de telles structures compatibles.

{\bf Compatibility betwenn pseudo-riemannian structure and Poisson structure}

{\bf Abstract. }We will introduce two notions of compatibility between  pseudo-riemannian metric and Poisson structure using the notion of contravariant connection introduced in [1], we will study some proprieties of manifold endowed with such compatible structures and we will give some examples.

\section{ Deux conditions naturelles de compatibilit\'e entre une structure pseudo-riemannienne et une structure de Poisson}
A notre connaissance, la notion de compatibilit\'e entre une structure pseudo-riemannienne et une structure de Poisson n'a jamais \'et\'e \'etudi\'e bien que ce soit une notion naturelle et bien qu'il  soit \'etabli que  la compatibilit\'e entre deux structures engendre des situations g\'eom\'etriques tr\`es riches ( compatibilit\'e entre deux structures de Poisson, compatibilit\'e entre une m\'etrique et une structure symplectique, vari\'et\'es K\"ahl\'eriennes). La difficult\'e pour trouver une telle notion est le fait que la m\'etrique et le tenseur de Poisson appartiennent \`a deux univers diff\'erents, l'un covariant et l'autre contravariant. Dans ce qui suit, nous allons donner deux conditions de compatibilit\'e ind\'ependantes en pr\'evilig\'eant le point de vue contravariant. Pour les notions de base sur les structures de Poisson voir, par exemple [4], et pour la notion de d\'eriv\'ee contravariante voir [1].
\subsection{D\'eriv\'ees contravariantes associ\'ees \`a une m\'etrique pseudo-riemannienne et un champ de bivecteurs}
Soit $M$ une vari\'et\'e diff\'erentiable munie d'une m\'etrique pseudo-riemannienne $g$ et d'un champ de bivecteurs $\pi$. 

Au champ de bivecteurs $\pi$ est associ\'e, d'une mani\`ere naturelle, un champ
d'endomorphismes $\#_{\pi}:T^*M\too TM$ par
$\be(\#_{\pi}(\al))=\pi(\al,\be)$
et  un crochet sur les 1-formes
d\'efini par
$$[\al,\beta]_\pi=L_{\#_\pi(\al)}\be-L_{\#_\pi(\be)}\al-d(\pi(\al,\be)).\eqno(1)$$
Le
crochet de Shouten-Nihenjuis $[\pi,\pi]_S$ est l'obstruction \`a ce que $\#_\pi$
soit un homorphisme entre $\Omega^1(M)$ munie du crochet $[,]_\pi$ et ${\cal
X}(M)$ muni de crochet de Lie. Nous avons
$$[\pi,\pi]_{S}(\al,\be,\ga)=\ga\left[\#_\pi([\al,\beta]_\pi)-
[\#_\pi(\al),\#_\pi(\beta)]\right].\eqno(2)$$		
Le champ de bivecteurs $\pi$ d\'efinit une structure de Poisson sur $M$ si et
seulement si
$[\pi,\pi]_S=0.$

 Soit $b_g:TM\too T^*M$ et  $\#_g:T^*M\too TM$ les isomorphismes fibr\'es d\'efinis par $g$,
 soit  $\wi g$ la m\'etrique sur le fibr\'e cotangent  d\'efinie par
$\wi g(\al,\be)=g(\#_g(\al),\#_g(\be))$ et soit  $J:T^*M\too T^*M$ le champ d'endomorphismes d\'efini par
$\pi(\al,\be)=\wi g(J\al,\be)=-\wi g(\al,J\be).$ On consid\`ere aussi le champ d'endomorphismes
$\wi J=\#_g\circ J\circ b_g.$

Maintenant, nous sommes en mesure d'associer au couple $(g,\pi)$ deux d\'eriv\'ees contravariantes. La premi\`ere, qu'on notera $\nabla^\pi$,  est d\'efinie \`a l'aide de la connexion de Levi-Civita $\nabla$ de $g$ par
$$\nabla^\pi_\al\be=\nabla_{\#_\pi(\al)}\be.\eqno(3)$$
La deuxi\`eme, qu'on notera
$D^\pi$, est donn\'ee par la formule
\begin{eqnarray*}
2\wi g(D^\pi_\al\be,\ga)&=&\#_\pi(\al).\wi g(\be,\ga)+\#_\pi(\be).\wi
g(\al,\ga)-\#_\pi(\ga).\wi g(\al,\be)\\
&&+\wi g([\ga,\al]_\pi,\be)+\wi g([\ga,\be]_\pi,\al)+\wi
g([\al,\be]_\pi,\ga).\qquad (4)\end{eqnarray*}
La proposition suivante donne les premi\`eres propri\'et\'es des deux connexions $\nabla^\pi$ et $D^\pi$ et permet de les comparer.
\begin{pr} Nous avons:	
			
$1)$ $\nabla^\pi$ et $D^\pi$ sont m\'etriques $ie$ 
 $$D^\pi \wi g(\al,\be,\ga)=\#_\pi(\al).\wi g(\be,\ga)-\wi g(D^\pi_\al\be,\ga)-\wi g(\be,D^\pi_\al\ga)=0$$ et la m\^eme formule pour $\nabla^\pi$;
			
$2)$ $D^\pi$ n'a pas de torsion,
par contre $\na^\pi$ a de la
torsion et nous avons
$$\wi g([\al,\be]_\pi-(\na_\al^\pi\be-\na_\be^\pi\al),\ga)=
g(\na_{\#_g(\ga)}(\wi J)(\#_g(\al)),\#_g(\be));$$
\begin{eqnarray*}
3)\quad[\pi,\pi]_S(\al,\be,\ga)&=&
\frac12\{\wi g(\al,D_\ga^\pi( J)(\be))+
\wi g(\be,D_\al^\pi( J)(\ga))
+\wi g(\ga,D_\be^\pi( J)(\al))\}\\
&=&\wi g(\al,\na_\ga^\pi( J)(\be))+
\wi g(\be,\na_\al^\pi( J)(\ga))+\wi g(\ga,\na_\be^\pi( J)(\al));\end{eqnarray*}

$4)$ On pose
$T(\al,\be)=[\al,\be]_\pi-(\na_\al^\pi\be-\na_\be^\pi\al)$ et 
$\wi{\na}^\pi_\al\be={\na}^\pi_\al\be+\frac12T(\al,\be).$
$\wi{\na}^\pi$ est une d\'eriv\'ee contravariante sans torsion et nous avons
$$\wi{\na}^\pi_\al\be-D^\pi_\al\be=\frac12b_g\left[\na_{\#_g(\al)}(\wi J)(\#_g(\be))+
\na_{\#_g(\be)}(\wi J)(\#_g(\al))\right];$$

$5)$ On pose $S(\al,\be)=\wi{\na}^\pi_\al\be-D^\pi_\al\be.$
Nous avons
$$\wi g(D^\pi_\al( J)(\be)-\frac12\na^\pi_\al(
J)(\be),\ga)+\frac12[\pi,\pi]_S(\al,\be,\ga)=\wi g(S(\al,\ga),\wi J\be)-
\wi g(S(\al,\be),\wi J\ga).$$
\end{pr}	
{\bf Preuve:} Le 1) et la premi\`ere partie de 2)  est un calcul  trivial, la deuxi\`eme partie de 2) est un calcul direct et le reste d\'ecoule facilement de ce cette formule.$\Box$
\subsection{D\'efinition  de la compatibilit\'e entre une structure pseudo-riemannienne et une structure de Poisson}
Les connexions $\na^\pi$ et $D^\pi$  sont m\'etriques, il est alors  naturel de poser la d\'efinition suivante:
\begin{Def} Soit $M$ une vari\'et\'e diff\'erentiable munie d'une m\'etrique pseudo-riemannienne et d'un champ de bivecteurs $\pi$.

1) On dira que le couple $(g,\pi)$ est $\na^\pi$-compatible si 
 $\na^\pi\pi=0$ ou, d'une mani\`ere \'equivalente, $\na^\pi(J)=0$.

2)  On dira que le couple $(g,\pi)$ est $D^\pi$-compatible si  $D^\pi\pi=0$
ou, d'une mani\`ere \'equivalente, $D^\pi(J)=0$.	
 \end{Def}

{\bf Remarque.} 1) Si le couple $(g,\pi)$ est compatible pour l'une ou l'autre et si $f$ est une fonction de Casimir de $\pi$ alors le couple $(g,f\pi)$ est compatible.

2) D'apr\`es les relations $3)$  de la Proposition 1.1, la condition de compatibilit\'e entraine l'int\'egrabilit\'e du tenseur $\pi$.

3) D'apr\`es 5) de la Proposition 1.1, il semblerait que ces deux notions soient ind\'ependantes. Nous allons voir sur un exemple \`a la fin de cette note qu'elles le sont effectivement.

4) 
 Supposons que  $\#_\pi$ soit inversible ce qui \'equivaut \`a $J$ ou $\wi J$ inversible.

 Le couple $(g,\pi)$ est $\na^\pi$-compatible
si  et seulement si la 2-forme $\om(u,v)=\pi(\#_\pi^{-1}(u),\#_\pi^{-1}(v))$ est parall\`ele pour la connexion de Levi-Civita de $g$.

En utilisant (4), il est facile d'\'etablir la formule
 $$D^\pi_\al\be=\#_\pi^{-1}\left(\na^{\wi J}_{\#_\pi(\al)}\#_\pi(\be)\right)\eqno(5)$$
o\`u $\na^{\wi J}$ est la connexion de Levi-Civita de la m\'etrique $g^{\wi J}(u,v)=g({\wi J}^{-1}(u),
{\wi J}^{-1}(v))$. Le couple $(g,\pi)$ est $D^\pi$-compatible
si et seulement si la 2-forme $\om(u,v)=\pi(\#_\pi^{-1}(u),\#_\pi^{-1}(v))$ est parall\`ele pour la connexion $\na^{\wi J}$ qui est aussi la connexion de Levi-Civita de $g$.

Ainsi les deux conditions coincident dans le cas symplectique et nous  retrouvons la notion de compatibilit\'e \'etudi\'ee dans [3].

\subsection{ Interpr\'etation g\'eom\'etrique de la compatibilit\'e du couple $(g,\pi)$}
	
Il est clair que le couple $(g,\pi)$ est $\na^\pi$-compatible si et seulement si le tenseur de Poisson est invariant  par transport parall\`ele ( associ\'e \`a la connexion de Levi-Civita de $g$) le long des courbes tangentes aux feuilles symplectiques.
	
A la d\'eriv\'ee contravariante $D^\pi$ est associ\'ee une notion de transport parall\`ele le long de  courbes cotangentes ( voir [1]). Une courbe cotangente est un couple $(\ga,\al)$ o\`u $\ga:[0,1]\too S$ est une courbe diff\'erentiable dans une feuille symplectique $S$ et $\al:[0,1]\too T^*M$ diff\'erentiable telle que, pour tout $t\in[0,1]$, $\#_\pi(\al(t))=\frac{d\ga}{dt}(t).$ A chaque courbe cotangente $(\ga,\al)$ est associ\'e un isomorphisme d'espaces vectoriels $\tau_{(\ga,\al)}:T^*_{\ga(0)}M\too 
T^*_{\ga(1)}M$.  Le  couple $(g,\pi)$ est $D^\pi$-compatible si et seulement si  $\tau_{(\ga,\al)}$ est une isom\'etrie qui commute avec $J$ et ce pour toute courbe cotangente et toute feuille symplectique.\bigskip
\section{Compatibilit\'e et feuilletage symplectique}
Soit $M$ une vari\'et\'e diff\'erentiable munie d'une m\'etrique $g$  et d'un champ de bivecteurs $\pi$. On reprend les notations et les d\'efinitions de la section pr\'ec\'edente.

On suppose que le couple $(g,\pi)$ est compatible  ( pour l'une ou pour l'autre des notions) et on consid\`ere $S$ une feuille symplectique du feuilletage d\'efini par la structure de Poisson et on note $\om_S$ la forme symplectique de $S$.

 Si le couple $(g,\pi)$ est $\na^\pi$-compatible, il est facile de voir que   $S$ est une sous-vari\'et\'e totalement g\'eod\'esique et la forme symplectique est parall\`ele pour la  restriction de $\na$  \`a $S$. 

On suppose maintenant que le couple $(g,\pi)$ est $D^\pi$-compatible. Pour tout couple de champs de vecteurs $X,Y$ tangents \`a $S$, on pose
$$\na^S_XY=\#_\pi(D^\pi_\al\be)\eqno(6)$$ o\`u
 $\al,\be$ sont deux 1-formes telles que $\#_\pi(\al)=X$ et $\#_\pi(\be)=Y$.
De la formule $\#_\pi(D^\pi_\al\be)=\#_\pi(D^\pi_\be\al)+[\#_\pi(\al),\#_\pi(\be)]$ et du fait que $D^\pi\pi=0$, on d\'eduit que $\#_\pi(\al)=0$ ou 
$\#_\pi(\be)=0$ entraine que $\#_\pi(D^\pi_\al\be)=0$. Il en d\'ecoule que $\na^S$ d\'efinit une connexion covariante sur $S$. Un calcul direct permet de montrer que $\na^S\om_S=0$ et si la restriction de $g$ \`a $S$ est non-d\'eg\'en\'er\'ee, $\na^S$ est la connection de Levi-Civita de la restriction de $g$ \`a $S$. Les feuilles symplectiques ne sont pas en g\'en\'eral  totalement g\'eod\'esiques.

{\bf Remarque.} Si le couple $(g,\pi)$ est compatible  ( pour l'une ou pour l'autre des notions) et si $J^3=-J$ alors toutes les feuilles symplectiques sont des
vari\'et\'es   K\"ahl\'eriennes.

\section{ Tenseur de courbure de $D^\pi$}
Soit $M$ une vari\'et\'e diff\'erentiable munie d'une m\'etrique $g$  et d'un champ de bivecteurs $\pi$. On reprend les notations et d\'efinitions de la Section 1.

La connexion contravariante $D^\pi$ est sans torsion et comme dans le cas covariant, on lui associe un tenseur de courbure qu'on notera $R^\pi$ et un tenseur de Ricci qu'on notera $r^\pi$ par les formules 
$R^\pi(\al,\be)=D^\pi_{[\al,\be]_\pi}-(D^\pi_\al D^\pi_\be-
D^\pi_\be D^\pi_\al)$ et $r^\pi(\al,\be)=Tr(\ga\too R^\pi(\al,\ga)\be).$ Ces deux tenseurs contravariants v\'erifient les m\^emes formules que dans le cas covariant.

  On suppose que le couple $(g,\pi)$ est $D^\pi$-compatible et on consid\`ere $S$ une feuille symplectique et on note $R^S$ et $r^S$ respectivement le tenseur de courbure et le tenseur de Ricci de la connexion $\na^S$ d\'efinie  par (6). Alors on a, clairement en vertu de (6)
$$R^S(X,Y)Z==\#_\pi(R^\pi(\al,\be)\ga)\eqno(7)$$
avec $X=\#_\pi(\al)$ et ainsi de suite.

Pour la courbure de Ricci, on remarque d'abord que $\#_\pi$ r\'ealise un isomorphisme entre $ImJ$ et $Im\wi J=TS$ et si on note 
$\#_\pi^{-1}$ sont isomorphisme inverse, un calcul direct nous donne la formule
$$
r^S(X,Y)=r^\pi(\#_\pi^{-1}(X),\#_\pi^{-1}(Y)).\eqno(8)$$

Ces deux formules sont int\'eressantes dans la mesure qu'elle permettent d'utiliser un tenseur global $R^\pi$ pour \'ecrire des propri\'et\'es des feuilles symplectiques.
 Par exemple, si $r^\pi=\lambda \wi g$ ( ici $\lambda$ est n\'ecessairement une fonction de Casimir de la structure de Poisson) alors toutes les feuilles symplectiques sont des vari\'et\'es D'Einstein.

Pour finir cette section, on remarque que, dans le cas de  $D^\pi$-compatibilit\'e, $r^\pi$ v\'erifie $r^\pi(\al,J\be)=-r^\pi(J\al,\be)$ et d\'efinit donc un champ de bivecteurs qu'on notera $\pi^r$ par la formule
$$\pi^r(\al,\be)=r^\pi(J\al,\be).\eqno(9)$$ Comme dans le cas K\"ahl\'erien, on a
$$[\pi,\pi^r]_S=0\eqno(10)$$ et on obtient donc une classe de cohomologie pour la cohomologie de Poisson.
\section{Exemples de couples compatibles}
Soit $(M,g)$ une vari\'et\'e pseudo-riemannienne et $U$ un champ de Killing de $g$. On d\'efinit le champ d'endomorphismes $\wi J:TM\too TM$ par
$$\wi JX=\na_XU\eqno(11)$$o\`u $\na$ est la connexion de Levi-Civita de $g$. $\wi J$ est anti-sym\'etrique par rapport \`a $g$ et si on consid\`ere $J=\#_g\circ \wi J\circ b_g$,
 on a un champ de bivecteurs par
$$\pi(\al,\be)=g(\wi J\#_g(\al),\#_g(\be)).$$On notera $R$ le tenseur de courbure de $g$.
\begin{pr} 1) Le couple $(g,\pi)$ est $\na^\pi$-compatible si et seulement si
$$R(U,\wi JX)Y=0\eqno(12)$$ pour tout couple de champs de vecteurs $X,Y$ sur $M$.

2)  Le couple $(g,\pi)$ est $D^\pi$-compatible si et seulement si
$$g(R(Y,Z)U,\wi JX)+g(R(X,U)Y,\wi JZ)+g(R(U,X)Z,\wi JY)=0\eqno(13)$$ pour tout triple de champs de vecteurs $X,Y,Z$ sur $M$.\end{pr}
{\bf Preuve:} 1) D\'ecoule imm\'ediatement du fait que $\na_X(\wi J)(Y)=R(U,X)Y$ pour tout couple de champs de vecteurs $X,Y$. 2) d\'ecoule aussi de cette relation et de la formule 5) Proposition 1.1.$\Box$

Un cas particulier o\`u on peut expliciter des exemples est le cas o\`u $M$ est un groupe de Lie $G$,  $g$ une pseudo- m\'etrique bi-invariante et $U$ un champ invariant \`a gauche associ\'e \`a un vecteur $u$ de l'alg\`ebre de Lie $\cal G$ de $G$. Dans ce cas, on a

\begin{pr} 1) Le couple $(g,\pi)$ est $\na^\pi$-compatible si et seulement si
$$[[u,[x,u]],y]=0\eqno(14)$$ pour tout couple  de vecteurs $x,y$ dans  $\cal G$.

2)  Le couple $(g,\pi)$ est $D^\pi$-compatible si et seulement si
$$[[u,x],[u,y]]=0\eqno(15)$$  pour tout couple  de vecteurs $x,y$ dans  $\cal G$.
\end{pr}
{\bf Preuve:} D\'ecoule de la Proposition 4.1 et du fait que dans ce cas $R(X,Y)=\frac14[X,Y]$ pour tout couple de champs de vecteurs $X,Y$ invariants \`a gauche.$\Box$

Les alg\`ebres de Lie munies de pseudo-m\'etriques bi-invariantes ont \'et\'e \'etudi\'ees dans [2]. Dans ce travail on trouve notamment l'exemple de $\reel^5 $ munie de la structure d'alg\`ebre de Lie donn\'ee par $[e_1,e_5]=e_4-e_3,$ $[e_2,e_5]=-e_1+e_3+e_4$,
$[e_1,e_2]=e_3$ et les autres tous nuls. Il est facile de voir que tout vecteur dans cette alg\`ebre de Lie v\'erifie (15).

D'un autre c\^ot\'e, il est facile de v\'erifier sur ce genre d'exemples que les conditions (14) et (15) sont ind\'ependantes.

{\bf R\'ef\'erences bibliograhiques}

[1] Fernandes R. L., Connections in Poisson Geometry1: Holonomy and invariants, J. of Diff. Geometry 54, 303-366, (2000).

[2] Medina A., Groupes de Lie munis de pseudo-m\'etriques de Riemann bi-invariantes, S\'eminaire Montpellier 1981-1982.

[3] Rakotondralambo J., Compatibilit\'e d'une forme symplectique et d'une pseudo-m\'etrique, th\`ese de l'Universit\'e de Pau, (D\'ecembre 1997).

[4] Vaisman I., Lecture on the geometry of Poisson manifold, Progr. In Math. Vol. 118, Birkhausser, Berlin, (1994).

\begin{tabular}{l}M. Boucetta\\
Facult\'e des Sciences et Techniques, BP 618, Marrakech, Maroc\\
Courriel: boucetta@fstg-marrakech.ac.ma\end{tabular}

\end{document}